\documentclass[11pt]{amsart}

\usepackage{lmodern,amssymb,microtype,fullpage,hyperref}
\usepackage[shortlabels]{enumitem}

\hypersetup{colorlinks=true, citecolor=blue, linkcolor=blue, urlcolor=blue, pdfstartview=FitH, pdfauthor=Gergely Harcos and Kannan Soundararajan, pdftitle=A supplement to Chebotarev's density theorem}

\newcommand{\CC}{\mathbb{C}}
\newcommand{\mo}{\mathfrak{o}}
\newcommand{\mpr}{\mathfrak{p}}

\newcommand{\ov}[1]{\overline{#1}}
\newcommand{\eps}{\varepsilon}

\renewcommand{\geq}{\geqslant}
\renewcommand{\leq}{\leqslant}
\renewcommand{\ge}{\geqslant}
\renewcommand{\le}{\leqslant}

\DeclareMathOperator{\Frob}{Frob}
\DeclareMathOperator{\Gal}{Gal}
\DeclareMathOperator{\Irr}{Irr}
\DeclareMathOperator*{\res}{res}
\DeclareMathOperator*{\ord}{ord}

\theoremstyle{plain}
\newtheorem*{maintheorem}{Theorem}
\newtheorem*{maincorollary}{Corollary}

\theoremstyle{remark}
\newtheorem{remark}{Remark}

\theoremstyle{definition}
\newtheorem*{acknowledgement}{Acknowledgements}

\begin{document}

\author{Gergely Harcos}
\address{Alfr\'ed R\'enyi Institute of Mathematics, POB 127, Budapest H-1364, Hungary}
\email{gharcos@renyi.hu}
\author{Kannan Soundararajan}
\address{Department of Mathematics, Stanford University, Stanford CA 94305, USA} 
\email{ksound@stanford.edu} 

\title{A supplement to Chebotarev's density theorem}

\dedicatory{On the 50th anniversary of Chen's theorem and the 100th anniversary of Chebotarev's theorem}

\begin{abstract} Let $L/K$ be a Galois extension of number fields with Galois group $G$. We show that if the density
of prime ideals in $K$ that split totally in $L$ tends to $1/|G|$ with a power saving error term, then the density
of prime ideals in $K$ whose Frobenius is a given conjugacy class $C\subset G$ tends to $|C|/|G|$ with the same power
saving error term. We deduce this by relating the poles of the corresponding Dirichlet series to the zeros of
$\zeta_L(s)/\zeta_K(s)$.
\end{abstract}



\thanks{The first author was supported by the R\'enyi Int\'ezet Lend\"ulet Automorphic Research Group and NKFIH (National Research, Development and Innovation Office) grant K~143876.  The second author was supported in part by a grant from the National Science Foundation, and a Simons Investigator Award from the Simons Foundation.}
\maketitle

\section{Introduction}

This note arose from an amusing observation by MathOverflow user Lucia~\cite{Mo}, which informally says that if the density of primes congruent to $1$ modulo $q$ tends to $1/\phi(q)$ rapidly, then for all $(a,q)=1$ the density of primes congruent to $a$ modulo $q$ also tends to $1/\phi(q)$ rapidly. Here density refers to the actual proportions. More precisely, for any $\sigma\ge 1/2$, the asymptotic
\begin{equation}\label{eq1}
\psi(x;q,1) = \psi(x)/\phi(q) + O_\eps(x^{\sigma+\eps})
\end{equation}
implies that 
\begin{equation}\label{eq2}
\psi(x;q,a) = \psi(x)/\phi(q) + O_\eps(x^{\sigma+\eps})\qquad\text{for all}\qquad (a,q)=1.
\end{equation}
The reason is simple. The relation \eqref{eq1} implies that the function
\[\sum_{\substack{\chi\bmod q\\\chi\neq\chi_0}}\frac{L^{\prime}}{L}(s,\chi)
=\int_1^{\infty}\bigl(\psi(x,\chi_0)-\phi(q)\psi(x;q,1)\bigr) \frac{s}{x^{s+1}}\,dx\]
extends analytically to the half-plane
\[H_\sigma:=\{s\in\CC:\Re(s)>\sigma\}.\]
That is, the product of the Dirichlet $L$-functions $L(s,\chi)$ ($\chi\neq\chi_0$) has no zero or pole in $H_\sigma$. As these $L$-functions are entire, none of them has a zero or pole in $H_\sigma$, and \eqref{eq2} follows easily.

Our goal is to show that this phenomenon persists in the context of Chebotarev's density theorem~\cite{Ch1,Ch2} (cf.\ \cite[Ch.~VII, Th.~13.4]{Ne}), even though the underlying Artin $L$-functions are only conjectured to be entire (Artin's conjecture). 

\begin{maincorollary} Let $L/K$ be a Galois extension of number fields, with $G$ denoting the Galois group.  Suppose $\sigma \ge 1/2$ is such that for any $\eps >0$ and $x\geq 2$ the asymptotic formula 
$$
\sum_{\substack{ N({\mpr}) \le x \\ \Frob(\mpr) =\{1 \} } }\log N(\mpr) = \frac{1}{|G|}  \sum_{N({\mpr}) \le x} \log N({\mpr}) + O( x^{\sigma+\eps}) 
$$ 
holds. Here the sums are over unramified  prime ideals of $\mo_K$, and $N(\mpr):=|\mo_K/\mpr|$ is the absolute norm of $\mpr$.  Then, the Artin $L$-functions $L(s, \chi)$ associated to 
the non-trivial irreducible characters of $G$ are analytic in the half-plane $H_\sigma$ and have no zeros there.   Consequently, for any conjugacy class $C$ of $G$ we have 
$$ 
\sum_{\substack{ N({\mpr}) \le x \\ \Frob(\mpr) =C } }\log N(\mpr) = \frac{|C|}{|G|} \sum_{N({\mpr}) \le x} \log N({\mpr}) + O( x^{\sigma+\eps}).
$$
\end{maincorollary} 

\begin{remark}
If $H\lhd G$ is a normal subgroup, then the above result applies equally well to the Galois extension $L^H/K$ in place of $L/K$. It follows that if the density of prime ideals in $K$ whose Frobenius lies in $H$ tends to $|H|/|G|$ with a power saving error term, then the density of prime ideals in $K$ whose Frobenius lies in $CH$ tends to $|CH|/|G|$ with the same power saving error term. We are grateful to the referee for pointing out this formal generalization.
\end{remark}

Lucia's observation pertains to the special case when $K={\mathbb Q}$ and $L = {\mathbb Q}(e^{2\pi i/q})$.  In the next section we shall formulate a more precise result 
from which the corollary above will follow.  As we shall see, the phenomenon admits a clean explanation by Heilbronn characters~\cite{He,FGM}, and leads to a transparent proof of the density theorem itself. In particular, we shall use a result of Foote--Murty~\cite[\S 3, Prop.]{FM} (cf.\ \cite[Prop.~2.1]{FGM}), which refines Aramata's theorem~\cite{Ara} (cf.\ \cite[Ch.~XVIII, Th.~8.4]{La}) and ultimately relies on Artin's reciprocity law~\cite{Art} (cf.\ \cite[Ch.~VI, Th.~5.5]{Ne}).

\begin{acknowledgement}
MathOverflow user Lucia is aware of this note and agrees that it be published in the current state and authorship. The authors are grateful to Lucia. The authors also thank Emmanuel Kowalski and the anonymous referees for valuable comments.
\end{acknowledgement}

\section{Statement of the result}

As above, let $L/K$ be a Galois extension of number fields with Galois group $G:=\Gal(L/K)$. For each conjugacy class $C\subset G$,  consider the Dirichlet series
\[
F(s,C):=\sum_{\Frob(\mpr)=C}\frac{\log N(\mpr)}{N(\mpr)^s}-\frac{|C|}{|G|}\sum_{\mpr}\frac{\log N(\mpr)}{N(\mpr)^s},
\]
which converges absolutely in the half-plane $H_1$.   Here $\mpr$ runs through the unramified prime ideals of $\mo_K$, and $N(\mpr):=|\mo_K/\mpr|$ is the absolute norm of $\mpr$.
The Chebotarev density theorem, in its weakest analytic form, states that $F(s,C)$ tends to a finite limit as $s\to 1+$.

For our purposes, $F(s,C)$ is essentially the same as
\[H(s,C):=-\frac{|C|}{|G|}\sum_{\substack{\chi\in\Irr(G)\\\chi\neq\chi_0}}\frac{L^{\prime}}{L}(s,\chi)\ov{\chi}(g_C),\qquad s\in H_1,\]
where $\chi$ runs through the non-trivial irreducible characters of $G$.  Here $L(s,\chi)$ is the Artin $L$-function of $\chi$, and $g_C$ is any element of the conjugacy class $C$. Indeed, the Schur orthogonality relation~\cite[Ch.~XVIII, Th.~5.5]{La} shows that $F(s,C)-H(s,C)$ is represented by an absolutely convergent Dirichlet series in $H_{1/2}$, hence 
it extends to a bounded analytic function in $H_\sigma$ for any $\sigma>1/2$. 

Recall Brauer's theorem~\cite[Th.~1]{Br} (cf.\ \cite[Ch.~XVIII, Th.~10.13]{La}), which shows that each $L(s,\chi)$ is meromorphic on $\CC$.  Hence $H(s,C)$ is meromorphic on $\CC$ with simple poles, and it follows that $F(s,C)$ extends meromorphically to $H_{1/2}$ with (at most) simple poles. A natural question arises: where are the poles and how large are the corresponding residues? The Chebotarev density theorem tells us that $s=1$ is not a pole.  Artin's conjecture and the Riemann hypothesis for Artin $L$-functions imply that $F(s,C)$ has no pole at all in the half-plane $H_{1/2}$. 

Our result relates the poles and residues of $F(s,C)$ in the half-plane $H_{1/2}$ to the zeros of $\zeta_L(s)/\zeta_K(s)$ in $H_{1/2}$. Since $L/K$ is Galois, 
Aramata's theorem shows that $\zeta_L(s)/\zeta_K(s)$ is an entire function.  Since $\zeta_K(s)$ and $\zeta_L(s)$ both have simple poles at $s=1$, we know moreover 
that $\zeta_L(s)/\zeta_K(s)$ is non-zero at $s=1$.  Thus the result stated below can be read as a supplement to (or refinement of) Chebotarev's density theorem.

\begin{maintheorem}
For any point $s_0\in H_{1/2}$, the following statements are equivalent:
\begin{enumerate}[(a)]
\item $s_0$ is a zero of $\zeta_L(s)/\zeta_K(s)$;
\item $s_0$ is a pole of $F(s,\{1\})$;
\item $s_0$ is a pole of $F(s,C)$ for some conjugacy class $C\subset G$;
\item $s_0$ is a zero or pole of $L(s,\chi)$ for some non-trivial $\chi\in\Irr(G)$;
\item $s_0$ is a zero of $L(s,\chi)$ for some non-trivial $\chi\in\Irr(G)$.
\end{enumerate}
Moreover,
\begin{equation} 
\label{4} 
\sum_C\frac{|G|}{|C|}\left|\res_{s=s_0} F(s,C)\right|^2=
\sum_{\substack{\chi\in\Irr(G)\\\chi\neq\chi_0}}\left(\ord_{s=s_0}L(s,\chi)\right)^2
\le\left(\ord_{s=s_0}\zeta_L(s)\right)^2-\left(\ord_{s=s_0}\zeta_K(s)\right)^2.
\end{equation}
\end{maintheorem}

\section{Proof of the theorem}

We begin by proving the key relation \eqref{4}. By our initial remarks,
\begin{equation}\label{5}
\res_{s=s_0} F(s,C)=\res_{s=s_0} H(s,C)=
-\frac{|C|}{|G|}\sum_{\substack{\chi\in\Irr(G)\\\chi\neq\chi_0}}\left(\ord_{s=s_0}L(s,\chi)\right)\ov{\chi}(g_C),
\end{equation}
hence
\begin{align*}
\sum_C\frac{|G|}{|C|}\left|\res_{s=s_0} F(s,C)\right|^2
&=\sum_C\frac{|C|}{|G|}\Biggl|\sum_{\substack{\chi\in\Irr(G)\\\chi\neq\chi_0}}
\left(\ord_{s=s_0}L(s,\chi)\right)\ov{\chi}(g_C)\Biggr|^2\\
&=\frac{1}{|G|}\sum_{g\in G}\Biggl|\sum_{\substack{\chi\in\Irr(G)\\\chi\neq\chi_0}}
\left(\ord_{s=s_0}L(s,\chi)\right)\ov{\chi}(g)\Biggr|^2.
\end{align*}
Now the first part of \eqref{4} follows from the Schur orthogonality relation~\cite[Ch.~XVIII, Th.~5.2]{La},
and the second part of \eqref{4} is a result of Foote--Murty~\cite[\S 3, Prop.]{FM} (cf.\ \cite[Prop.~2.1]{FGM}) upon noting that $\zeta_K(s)=L(s,\chi_0)$. See also Remark~\ref{remark2} below.

The rest of the Theorem is straightforward. The implication $(a)\Rightarrow (b)$ follows from the factorization
\begin{equation}\label{6}
\zeta_L(s)/\zeta_K(s)=\prod_{\substack{\chi\in\Irr(G)\\\chi\neq\chi_0}}L(s,\chi)^{\chi(1)}.
\end{equation}
Indeed, if $s_0$ is a zero of this meromorphic function, then it is a pole of the logarithmic derivative of the function,
which is $-|G|H(s,\{1\})$. Hence $s_0$ is a pole of $H(s,\{1\})$, and of $F(s,\{1\})$ as well.

The implication $(b)\Rightarrow (c)$ is trivial. 
 
The identity \eqref{5}, or alternatively the first part of \eqref{4}, ensures that if some $F(s,C)$ has a pole at $s_0$, then $\ord_{s=s_0} L(s,\chi)$ 
must be non-zero for some non-trivial $\chi\in\Irr(G)$.  Thus $(c)$ implies $(d)$. 
 
If $L(s,\chi)$ has a zero or pole at $s_0$ for some non-trivial $\chi\in\Irr(G)$, then the second part of \eqref{4} shows that the right-hand side of \eqref{4} is strictly positive. Hence $s_0$ is a zero of $\zeta_L(s)/\zeta_K(s)$, and by \eqref{6} it is a zero of $L(s,\chi)$ for some non-trivial $\chi\in\Irr(G)$ as well. Thus $(d)$ implies $(e)$, and the same argument also gives that $(e)$ implies $(a)$.    
  
The proof of the Theorem is complete.

\begin{remark} There are other natural ways to see the equivalence of the five statements in the Theorem. For example, \eqref{4} can be perceived as the analytic embodiment of $(c)\Leftrightarrow(d)\Rightarrow(a)$. We have seen that $(a)\Rightarrow(b)$ follows from \eqref{6}, while $(b)\Rightarrow(e)$ clearly follows from \eqref{5}, and $(e)\Rightarrow(d)$ is trivial.
\end{remark}

\begin{remark}\label{remark2} The crucial input in the proof above is the result of Foote and Murty that we quote. This result can be formulated as the bound $\langle\theta_G,\theta_G\rangle_G\le\theta_G(1)^2$, where
\[\theta_G:=\sum_{\chi\in\Irr(G)}\left(\ord_{s=s_0}L(s,\chi)\right)\chi\]
is the Heilbronn character of $G$ relative to $s_0$. Foote and Murty assume that $s_0\neq 1$, and they actually prove that $|\theta_G(g)|\le\theta_G(1)$ holds for all $g\in G$. This stronger bound relies on two fundamental properties. First, the restriction of $\theta_G$ to a subgroup $H\le G$ equals $\theta_H$. Second, the Heilbronn character of a cyclic group is a character (not just a virtual character) by Artin's reciprocity law. Indeed, putting $H=\langle g\rangle$, and using the first property, then the second property, and then again the first property, we obtain $|\theta_G(g)|=|\theta_H(g)|\le\theta_H(1)=\theta_G(1)$. For $s_0=1$ the same argument works if we replace $\theta_G$ by
\[\theta_G+\chi_0=\sum_{\substack{\chi\in\Irr(G)\\\chi\neq\chi_0}}\left(\ord_{s=1}L(s,\chi)\right)\chi.\]
The conclusion is that this modified Heilbronn character is zero, which is equivalent to Chebotarev's density theorem. Of course Artin's reciprocity law is also present in our overall discussion via Brauer's theorem that Artin $L$-functions are meromorphic.
\end{remark}

\begin{remark}
One can also apply Heilbronn's idea to the familiar function
\[U_G(s,g):=\sum_{\chi\in\Irr(G)}\frac{L^{\prime}}{L}(s,\chi)\ov{\chi}(g),\qquad s\in H_1,\quad g\in G.\]
Namely, by Frobenius reciprocity, the restriction of $U_G(s,\ast)$ to a subgroup $H\leq G$ equals $U_H(s,\ast)$.
Hence, by Artin reciprocity, $U_G(s,g)$ is meromorphic on $\CC$ with simple poles, and it has a simple pole at $s=1$ with residue $-1$.
Using also
\[\frac{L^{\prime}}{L}(s,\chi)=\langle U_G(s,\ast),\ov{\chi}\rangle,\]
one obtains a weak form of Brauer's theorem: each $L^{\prime}(s,\chi)/L(s,\chi)$ is meromorphic on $\CC$ with simple poles, and for $\chi\neq\chi_0$, the point $s=1$ is not a pole. From this, one can derive a slightly weaker version of the main Theorem
without using Brauer's theorem: $(a)$, $(b)$, $(c)$ are equivalent to each other and to the statement that $s_0$ is a pole of $L^{\prime}(s,\chi)/L(s,\chi)$ for some non-trivial $\chi\in\Irr(G)$. For a self-contained argument along these lines, see the lecture slides at \href{https://ntrg.math.unideb.hu/GH2023Talk.pdf}{https://ntrg.math.unideb.hu/GH2023Talk.pdf}.
\end{remark}

\section{Deducing the Corollary}
Note that 
\[F(s,\{1\}) = \int_{1}^{\infty} \Biggl( \sum_{\substack { N(\mpr) \le x\\ \Frob(\mpr) = \{ 1\} }} \log N(\mpr) - \frac{1}{|G|} \sum_{N(\mpr) \le x} \log N(\mpr)\Biggr) \frac{s}{x^{s+1}}\,dx,\]
so that the hypothesis in the Corollary implies that $F(s,\{1\})$ extends analytically to the half-plane $H_\sigma$. By the implication $(d)\Rightarrow(b)$ of our Theorem, it follows that 
all the Artin $L$-functions $L(s,\chi)$ with non-trivial $\chi\in\Irr(G)$ extend analytically to this half-plane $H_\sigma$ and have no zeros there.  Arguing as in the 
proof of the prime number theorem, it follows that for every non-trivial $\chi\in\Irr(G)$ one has 
\[\sum_{N(\mpr) \le x} \chi(\Frob(\mpr)) \log N(\mpr) = O(x^{\sigma+\eps})\]
for any $\eps >0$. The Schur orthogonality relation~\cite[Ch.~XVIII, Th.~5.5]{La} gives 
$$ 
\sum_{\substack{N(\mpr) \le x\\ \Frob(\mpr)=C}} \log N(\mpr) - \frac{|C|}{|G|} \sum_{N(\mpr)\le x} \log N(\mpr)  = 
\frac{|C|}{|G|} \sum_{\substack{\chi\in\Irr(G)\\\chi\neq\chi_0}} \overline{\chi}(g_C) \sum_{N(\mpr)\le x} \chi(\Frob(\mpr)) \log N(\mpr), 
$$ 
and the conclusion of the Corollary follows.

\end{document}